\documentstyle[12pt]{article}

\newenvironment{Pf}{\par\noindent{\sc
Proof.~}}{\ $\Box$\par\smallskip}

\newtheorem{defn}{Definition}[section]
\newtheorem{lem}[defn]{Lemma}
\newtheorem{thm}[defn]{Theorem}
\newtheorem{prop}[defn]{Proposition}
\newtheorem{cor}[defn]{Corollary}
\newtheorem{rem}[defn]{Remark}

\newtheorem{assu}[defn]{Assumption}
\newtheorem{constr}[defn]{Construction}
\newtheorem{nota}[defn]{Notation}
\makeatletter
\def\theequation{\thedefn.\@arabic\c@equation}
\makeatother

\newcommand{\Bl}{B\ell}
\newcommand{\Sym}{\Sigma}

\def\re^#1_#2{a^{#1}_{#2}}
\newcommand{\al}{\alpha}
\def\reb^#1_#2{(n_{#2}a^{#1}_{#2})/m_{#1}}
\newcommand{\iso}{\cong}
\newcommand{\N}{{\bf N}}
\newcommand{\C}{{\bf C}}
\newcommand{\Z}{{\bf Z}}

\newcommand{\ichi}{{i,\chi}}
\renewcommand{\P}{{I\!\! P}}
\renewcommand{\O}{{\cal O}}
\newcommand{\I}{{\cal I}}
\renewcommand{\L}{{\cal L}}
\newcommand{\Y}{{\cal Y}}
\newcommand{\M}{{\cal M}}
\newcommand{\X}{{\cal X}}
\newcommand{\F}{{\cal F}}
\newcommand{\hilb}{\hbox{\rm Hilb}}
\newcommand{\Hilb}{\hbox{\rm Hilb}}
\newcommand{\longra}{\longrightarrow}

\newcommand{\NS}{N\!S}
\def\xto.#1{x_1,\ldots,x_{#1}}
\def\yto.#1{y_1,\ldots,y_{#1}}

\newcommand{\barx}{{X'}}
\newcommand{\barc}{{C'}}
\newcommand{\barf}{{F'}}

\title{On the Hilbert scheme of curves in higher-dimensional projective space}
\author{Barbara \hbox{\rm Fantechi${}^*$} \quad --- \quad Rita
Pardini\thanks{Both authors are members of GNSAGA of CNR.}}
\date{}
\begin{document}
\maketitle
\begin{abstract}\noindent In this paper we prove that, for any
$n\ge 3$, there exist infinitely many $r\in \N$ and for each of them a smooth,
connected curve
$C_r$ in
$\P^r$ such that
$C_r$ lies on  exactly $n$ irreducible components of the Hilbert scheme
$\hilb(\P^r)$. This is proven by reducing the problem to an analogous
statement for the moduli of surfaces of general type.
\end{abstract}
\section{Introduction}
It is well-known that the Hilbert scheme parametrizing
subschemes of
$\P^r$ can be singular at points corresponding to smooth curves as soon as
$r\ge 3$; actually Mumford \cite{Mu} gave an example of an everywhere singular
irreducible component. If $r=3$, it has been proven in \cite{EHM} that the
open subset of the Hilbert scheme parametrizing smooth curves in $\P^3$ with
given genus and degree can have arbitrarily many components when the genus and
the degree grow (in fact, they prove that no polynomial estimate on the number
of such components holds).

Our main result is the following:

\smallskip
\noindent
{\bf Theorem \ref{mainthm2}.}{\em\ Let $n\ge 3$ be an integer. Then there
exist infinitely many integers $r$, and for each of them a smooth, irreducible
curve
$C_r\subset \P^r$ such that $C_r$ lies exactly on $n$ components of the
Hilbert scheme of $\P^r$.}

\smallskip
The idea of the proof is very simple. Firstly, we modify a construction of
\cite{FP} to obtain a regular surface $S$ of general type which lies on $n$
components of the moduli space; secondly, we
consider a suitable pluricanonical embedding of this surface and intersect its
image with a high-degree hypersurface $F$ to construct the curve $C$ we are
interested in. Finally, we prove that all embedded deformations of $C$
are induced by embedded deformations of $F$ and $S$.

\par\noindent
{\em Acknowledgements}. We are grateful to Ciro
Ciliberto, who told us about this problem and suggested that we might apply
to it the results of \cite{FP}.
\section{Notation and preliminaries}
All varieties will be assumed smooth and projective over the complex numbers
unless the contrary is explicitly stated. A variety $Y$ will be called
regular if $H^1(Y,\O_Y)=0$. If $\F$ is a sheaf on $Y$, let $h^i(Y,\F)=\dim
H^i(Y,\F)$. If $t$ is a real number, we denote its integral part by $[t]$.
Let $\zeta_3=\exp(2\pi i/3)$.

In this paper we will be concerned with abelian covers of a very special
type; we collect here the necessary notational set-up.

Let $n$ be an integer $\ge 2$, and let $G=\Z_3^n$, $G^*$ its dual; let
$e_1,\ldots,e_n$ be the canonical basis of $G$, and $\chi_1,\ldots,\chi_n$
the dual basis of $G^*$ (i.e., $\chi_j(e_i)=1$ if $i\ne j$ and
$\chi_j(e_j)=\zeta_3$). Let $e_0=-(e_1+\ldots+e_n)$. Let
$I=\{0,\ldots,n\}$,
and to each
$i\in I$ associate the pair
$(H_i,\psi_i)$ where
$H_i$ is the cyclic subgroup of $G$ generated by $e_i$, and $\psi_i\in
H_i^*$ is the character such that $\psi_i(e_i)=\zeta_3$.

Let
$Y$ be a smooth projective variety, and $(G,I)$ as above: a $(G,I)$-cover of
$Y$ is a normal variety $X$ and a Galois cover $f:X\to Y$ with
Galois group $G$ and (nonempty) branch divisors $D_i$ (for $i\in I$) having
$(H_i,\psi_i)$ as inertia group and induced character (see \cite{Pa} for
details).
\begin{lem}
To give a smooth $(G,I)$-cover of $Y$ is equivalent to giving
line bundles $L$ and $F_j$, for $j=1,\ldots,n$, together with smooth
nonempty divisors $D_i\in |M_i|$ {\rm(}where $M_0=L$ and, for $i\ge 1$,
$M_i=L-3F_i${\rm)} such that the union of the $D_i$'s has normal crossings.
\end{lem}
\begin{Pf}
{}From \cite{Pa} we know that the cover is determined by its reduced building
data, divisors $D_i$ for $i\in I$ and line bundles $L_j$ for $j=1,\ldots,n$
satisfying the relation  $3L_j\equiv D_j+2D_0$. Letting $M_i=\O(D_i)$, and
putting $M_0=L$, $F_j=L-L_j$, the equations become precisely $M_j=L-3F_j$.
\end{Pf}

As the natural map
$\bigoplus_{i\in I}H_i\to G$ is
surjective, the covers we consider will be totally ramified. For $\chi\in
G^*$, let as usual $L_\chi^{-1}$ be the corresponding eigensheaf in the
direct sum decomposition of $f_*\O_X$; in the above notation, we will have
(for $\chi=\chi_1^{\al_1}\cdots\chi_n^{\al_n}$):
\setcounter{equation}{0}
\begin{equation}\label{lchi}
L_\chi=n_\chi L-\sum_{j=1}^n \al_jF_j,
\end{equation}
where $n_\chi=-[(-\al_1-\ldots-\al_n)/3]$. In particular note
that $n_\chi\ge 1$ when $\chi\ne 1$, and $n_\chi=1$ if and only if $1\le\sum
\al_j\le 3$. We will write $L_j$ instead of $L_{\chi_j}$.

Recall from \cite{Pa}, proof of proposition 4.2 on page 208, that
\begin{equation}\label{canonico}
3K_X=\pi^*(3K_Y+2(n+1)L-6\sum F_j).
\end{equation}
\smallskip
We now recall some results from \cite{FP} in a simplified form (fit for our
situation). For details and proofs see \cite{FP}, \S 5.

\begin{rem}\label{rema}{\rm (1) Let $\Y\to B$ be a
smooth projective morphism (with $B$ a smooth, connected quasiprojective
variety) together with an isomorphism between $\Y_o$ and $Y$ for some $o\in
B$, and assume that $Y$ is regular and that the morphism $\Y\to B$ has a
section $\sigma$. Let
$L$ be a line bundle on
$Y$; assume that $c_1(L)$ is kept fixed by the
monodromy action of
$\pi_1(B,o)$ on
$H^2(Y,\Z)$. Then for each
$b\in B$ there is a canonical induced class
$c_1(L_b)$ on $\Y_b$. If, for all $b\in B$, the class $c_1(L_b)$ is of type
$(1,1)$, then
$L$ can be extended to a line bundle $\L$ over $\Y$, flat over $B$; this
extension is unique if we require that its restriction to $\sigma(B)$ be
trivial. This follows by applying the results on p.~20 of \cite{mumford}, and
by noting that the relative Picard scheme of $\Y$ over $B$
is \'etale over $B$ since all fibres are smooth and regular (it is surjective
as $c_1(L_b)$ is always of type $(1,1)$); the condition on the monodromy action
implies then that the component of the relative Picard scheme containing $[L]$
is in fact isomorphic to $B$. Let
$L_b$ be the restriction of
$\L$ to
$\Y_b$.

\noindent (2) If $h^0(\Y_b,L_b)$ is either constant in $b$, or if it only
assumes the values $1$ (for $b\in Z$) and $0$, then there is a
(nonunique) quasiprojective variety $W^L\to B$ such that $W^L_b$ is canonically
isomorphic to $H^0(\Y_b,L_b)$; $W^L$ is smooth and irreducible in the former
case, while in the latter it is the union of one component isomorphic to $B$
and another being the total space of a line bundle over $Z$ (compare with
\cite{FP}, theorem 5.8 and remark 5.11). }\end{rem}

\begin{assu}\label{ass}{\rm
Let $S=\{(\ichi)\in I\times G^*|\chi_{|H_i}\ne\psi_i^{-1}\}$.
Let $X\to Y$ be a smooth $(G,I)$-cover as in lemma 2.1, and
$\Y\to B$ be a smooth projective morphism (with $(B,o)$ a pointed space, and
$\Y_o$ isomorphic to $Y$), such that remark
\ref{rema}, (1) applies to $\Y\to B$, for the line bundle $L$ and for each of
the $F_j$'s.  Assume moreover that remark \ref{rema}, (2) applies for the line
bundles
$M_i-L_\chi$ for $(\ichi)\in S$, yielding varieties $W^\ichi$: let $W$ be the
fibred product of the $W^\ichi$ over $B$. Finally, assume that the germ of $B$
at $o$ maps smoothly to the base of the Kuranishi family of $Y$, and that the
cohomology groups $H^1(Y,L_\chi^{-1})$ and $H^1(Y,T_Y\otimes L_\chi^{-1})$
vanish
for each $\chi\in G^*\setminus 1$.
}\end{assu}

\begin{thm}\label{fromFP}
Assume that assumption {\rm\ref{ass}} holds, and let
$w\in W$ be a point over $o\in B$ corresponding to sections $s_\ichi$ such that
$s_\ichi=0$ if $\chi\ne 1$, and $s_{i,1}$ defines $D_i$ for $i=0,\ldots,n$.
Assume also that $X$ has ample canonical class. One can construct a family of
natural deformations of
$(G,I)$-covers
$\X\to W$; the induced map from the germ of $w$ in $W$ to the Kuranishi family
of
$X$ is smooth {\rm(}and, in particular, surjective{\rm)}. Moreover, the flat,
projective morphism
$\X\to W$ defines a rational map from $W$ to the moduli of surfaces with
ample canonical class, regular at $w$; this map is
dominant on each irreducible component of the moduli containing $X$.
\end{thm}

\begin{Pf}
Let $\L$ (resp.~$\F_j$) be the line bundle induced by $L$
(resp.~$F_j$) on $W$; as we {\em define} $\M_0$ to be $\L$, $\L_j$ to be
$\L-\F_j$ and  $\M_j$ to be $\L-3\F_j$ for $j=1,\ldots,n$, there are global,
canonical isomorphisms $\phi_j:3\L_j\to \M_j+2\M_0$.
By \cite{FP}, theorem 5.12, the germ of $W$ at $w$ maps smoothly to the
base of the Kuranishi family of $X$.

If $M$ is
an irreducible component of the moduli containing $[X]$, by the previous result
the image of $W$ contains an open set in $M$ (in the strong topology), hence it
cannot be contained in a closed subset (in the Zariski topology) and is
therefore dominant.
\end{Pf}

\section{Moduli of surfaces of general type}
The aim of this section is the proof of theorem \ref{mainthm1}, i.e., the
explicit construction of regular surfaces with ample canonical class lying on
arbitrarily many components of the moduli. This construction can be carried
out in a much more general setting (see remark \ref{+gen}); we consider only
the case needed for our applications, since it is easier to describe.

For $S$  a smooth projective surface and
$\xto.n$  pairwise distinct points of $S$, we let $\Bl(S;\xto.n)$ denote the
surface obtained by blowing up $S$ at $\xto.n$.
\begin{constr}{\rm Let $S$ be a regular surface, $x_0\in S$, $n$ a positive
integer; let
$B=B(S,n)$ be the variety parametrizing data $(\xto.n,\yto.n)$ where
the $x_i$'s are pairwise distinct points in $S$ (for $i=0,\ldots,n$), the
$y_i$'s
are pairwise distinct points in $\Bl(S;\xto.n)$, such that $y_i$ is not
infinitely near to $x_j$ for $i\ne j\ge 1$ and none of the $y_i$'s lies
over $x_0$.
$B$ is a smooth quasiprojective variety, which is naturally isomorphic to an
open subset of the product of $n$ copies of
$S\times S$ blown up along the diagonal. Let $\Y\to B$ be the smooth
projective family such that $\Y_b$, the fibre of $\Y$ over the point $b$, is
isomorphic to
$\Bl(\Bl(S;\xto.n);\yto.n)$ for $b=(\xto.n,\yto.n)$.}
\end{constr}
Note that the morphism $\Y\to B$ has a section, given by mapping $b\in B$ to
the inverse image of $x_0$ in $\Y_b$.
\begin{lem} Assume that
$S$ is rigid. Let $B^0$ be the open set in $B$ where
$Aut(S)$ acts freely (the action being the natural one). Then if $b\in B^0$,
the natural map from the germ of $B$ in $b$ to the Kuranishi family of $\Y_b$
is smooth of relative dimension $h^0(S,T_S)$.
\end{lem}
\begin{Pf} The proof is easy and left to the reader.
\end{Pf}

\begin{rem}{\rm For any $b\in B$, $b=(\xto.n,\yto.n)$, there is a canonical
isomorphism
$$\NS(\Y_b)=\NS(S)\oplus\Z e'_1\oplus\ldots\oplus\Z e_n'\oplus\Z e''_1\oplus
\ldots\oplus\Z e_n'',$$
where $e_i'$ is the pullback from $\Bl(S;\xto.n)$ of the class of the
exceptional divisor over $x_i$, and $e_i''$ is the class of the exceptional
divisor over $y_i$. We will consider this isomorphism fixed, and denote this
group by $\NS$. We also let $f_i$ denote $e_i'-e_i''$.
Since $S$ is regular, so are all the $\Y_b$'s and
we will not need to distinguish between line bundles and their Chern
classes.
}\end{rem}

\begin{defn}{\rm Let $L\in \NS$, $G=\Z_3^n$ as in \S
2; for $\chi\in G^*$, let $L_\chi\in \NS$ be defined by equation (\ref{lchi}),
with $F_i=f_i$. Let $B_{L}$ be the open subset
of $B$ consisting of the $b$'s such that
\begin{enumerate} \item the cohomology groups $H^1(\Y_b,L_\chi^{-1})$,
$H^1(\Y_b,T_{\Y_b}\otimes L_\chi^{-1})$ are zero for each $\chi\in G^*\setminus
1$;
\item the line bundles $L$ and $L-3F_j$ are very ample on $\Y_b$,
for $j=1,\ldots,n$;
\item the line bundles $L-K_{\Y_b}$ and $L-3F_j-K_{\Y_b}$ are ample on  $\Y_b$,
for $j=1,\ldots,n$;
\item the line bundle $3K_Y+2(n+1)L-6\sum F_j$ is ample on $\Y_b$.
\end{enumerate}}
\end{defn}
Note that the first condition is needed to ensure that assumption \ref{ass}
is satisfied; the second allows one to choose smooth divisors in the linear
systems $|L|$ and $|L-3F_j|$ meeting transversally; the third implies that
these linear systems have constant dimension when $b$ varies; and
the fourth ensures, in view of equation (\ref{canonico}), that the cover so
obtained has ample canonical class (recall that the pullback of an ample line
bundle via a finite map is again ample).

\begin{lem}\label{basic}
Let $Y$ be a smooth surface containing $m$
disjoint irreducible curves $C_1,\ldots,C_m$, such that $C_i^2<0$. Then:
\begin{enumerate}
\item for any choice of nonnegative integers $a_1,\ldots,a_m$,
the linear system $|a_1C_1+\ldots+a_mC_m|$ contains only the divisor
$a_1C_1+\ldots+a_mC_m$;
\item for any choice of nonnegative integers $a_1,\ldots,a_{m-1}$, and for
any $b>0$, the linear system $|a_1C_1+\ldots+a_{m-1}C_{m-1}-bC_m|$ is empty.
\end{enumerate}\end{lem}
\begin{Pf} (1) We prove the theorem by induction on $a_1+\ldots+a_m$, the
case where this sum is zero being trivial. Assume without loss of generality
that $a_1\ge 1$, and let $C\in  |a_1C_1+\ldots+a_mC_m|$; then $C\cdot
C_1=a_1C_1^2<0$, hence $C$ must have a common component with $C_1$; therefore
$C=C_1+C'$, $C'\in |(a_1-1)C_1+\ldots+a_mC_m|$, and by induction the proof is
complete.\par
\noindent
(2) Assume that there exists $C\in |a_1C_1+\ldots+a_{m-1}C_{m-1}-bC_m|$. Then
$C+bC_m\in |a_1C_1+\ldots+a_{m-1}C_{m-1}|$, contradicting (1).
\end{Pf}

\begin{cor} Let $b\in B$, $b=(\xto.n,\yto.n)$ and let $a_1,\ldots,a_m$ be
nonnegative integers. Then the line bundle
$a_1f_1+\ldots+a_mf_m$ on $\Y_b$ is effective if and only if $y_i$ is
infinitely near to
$x_i$ for every $i$ such that $a_i>0$, and in this case it has only one
section.
\end{cor}
\begin{Pf} If $y_i$ is infinitely near to $x_i$, then $f_i$ is a $(-2)$
curve; otherwise it is the difference of two disjoint $(-1)$ curves. In the
former case lemma \ref{basic}, (1) applies and in the latter case
\ref{basic}, (2) applies.
\end{Pf}
\begin{nota}{\rm We will denote by $E$ the closed subset of $B$
consisting of the points $b$ such that $f_i$ is effective on $\Y_b$ for
$i=1,\ldots,n$.}
\end{nota}

\begin{lem} Let $L\in NS$ be a line bundle and assume that $E\cap
B_L\ne\emptyset$. Then assumption {\rm\ref{ass}} holds for the restriction of
$\Y\to B$ to $B_L$; applying theorem {\rm\ref{fromFP}} yields a quasiprojective
variety $W$. In this case
$W$ is the union of
$2^n$ smooth irreducible components
$W_A$, indexed by subsets
$A\subset
\{1,\ldots,n\}$. The dimension of $W_A$ and $W_{A'}$ are equal if $\#A=\#A'$,
 and in particular one has:
$$\dim W_A-\dim W_{\emptyset}=\frac{1}{6}(\#A^3+6\#A^2-\#A).$$
The $W_A$'s have a nonempty intersection.
\end{lem}
\begin{Pf}
The verification that assumption \ref{ass} holds is easy and we leave it to
the reader. For $A\subset \{1,\ldots,n\}$, let $E_A=\{b\in B_L|f_i\ \mbox{is
effective for}\ i\in A\}$ and let $W_A\subset W$ be defined by
$$W_A=\{(b,s_{\ichi})|b\in E_A\ \mbox{and}\ s_{\ichi}=0\ \mbox{for}\ \chi \ne
1\ \mbox{and}\ i\notin A\}.$$ It is easy to check that $W_A$ is smooth over
$E_A$ of dimension $1/6(\#A^3+6\#A^2+5\#A)$; on the other hand $E_A$ is
smooth of codimension $\#A$ in $B_L$. Finally, $W$ is the union of the
$W_A$'s, which are easily seen to be irreducible components.
Moreover, the intersection of the $W_A$'s is clearly equal to
$W_\emptyset$ intersected with the inverse image of $E\cap B_L$.
\end{Pf}

\begin{thm}\label{mainthm1}
Let $L\in \NS$ and $b\in E\cap B_L\cap B^0$. Let
$f:X\to Y=\Y_b$ be a smooth
$(G,I)$-cover with building data $(D_i,L_\chi)$; let $w\in W_b$ be a point
corresponding to a choice of equations $s_i\in \O_Y(D_i)$ defining $D_i$,
with $s_{\ichi}=0$ for all $\chi\ne 1$. Then the natural map from the germ of
$W$ in $w$ to the Kuranishi family of $X$ is smooth. In particular, the
Kuranishi family of $X$ is the union of $2^n$ irreducible components,
$n+1$ of which have pairwise different dimension.
Moreover, the surface $X$ lies on exactly $n+1$ components of the
moduli space, having pairwise different dimensions.
\end{thm}
\begin{Pf}
The first statement is a straightforward application of theorem \ref{fromFP},
in view of the previous lemma.
To prove the second, note that the map from $W$ to the moduli space factors
through the action of the symmetric group on $n$ letters, $\Sym_n$. The
quotient $W/\Sym_n$ has exactly $n+1$ irreducible components of pairwise
different dimensions. By the previous result, each of these components
dominates a component of the moduli; the $n+1$ components so obtained must
all be distinct, as they have different dimensions.
\end{Pf}
\begin{rem}\label{suffample}
{\rm
If $L\in NS$ is sufficiently ample, then the intersection $E\cap B_L\cap B^0$
is nonempty, hence the theorem applies yielding infinitely many surfaces with
different invariants. }
\end{rem}
\begin{cor}\label{HilbX}
Given integers $n\ge 2$ and $m\ge 5$, for infinitely
many values of $r$ there exists a smooth, regular surface $X$ in $\P^r$
such that $\O_X(1)=mK_X$ and $X$ lies on exactly $n+1$ irreducible components
of the Hilbert scheme.
\end{cor}
\begin{Pf}
Let $X$ be a regular surface,
with ample $K_X$, lying on exactly $n+1$ irreducible components of the moduli;
let $M$ be the union of the irreducible components of the moduli space of
surfaces with ample canonical class which contain $[X]$. Let
$r=h^0(X,mK_X)-1$; infinitely many such $X$'s (with distinct values of $r$) can
be constructed by applying theorem \ref{mainthm1}, in view of remark
\ref{suffample}.

Fix an $m$-canonical embedding of $X$ in $\P^r$. Every small embedded
deformation of $X$ in $\P^r$ is again a smooth surface, $m$-canonically
embedded
as
$X$ is regular. Let $H$ be the union of the irreducible components of the
Hilbert
scheme of $\P^r$ containing $[X]$, and $H^0$ be the open dense subset of $H$
parametrizing smooth, $m$-canonically embedded surfaces.

The natural map $H^0\to M$ is dominant, and each fibre is irreducible of
dimension $(r+1)^2-1$; in fact, the fibre over $[X']$ is the set of
bases of $H^0(X',mK_{X'})$ modulo the action of the finite group $Aut(X')$ (and
modulo the obvious $\C^*$-action).

In particular there is an induced bijection between irreducible
components of $M$ and of $H$, which increases the dimension by $(r+1)^2-1$.
\end{Pf}
\begin{rem}\label{+gen}
{\rm The constructions in this section generalize easily to the case where $Y$
is neither regular nor rigid and $G$ is any abelian group. In fact, they also
work if the dimension of $Y$ is bigger than $2$ (using a suitable, modified
form of lemma \ref{basic}).}
\end{rem}

\section{The Hilbert scheme of curves in $\P^r$}
In this section we apply the results on the Hilbert scheme of surfaces to the
Hilbert scheme of curves. We first introduce some notation. If $Z\subset
\P^r$ is a subscheme, we will denote by $\Hilb(Z)$ the union of the
irreducible components of the Hilbert scheme of $\P^r$ containing $[Z]$; we
let $H(Z)$ be the germ of $\Hilb(Z)$ at $[Z]$.
Note that if $F$ is a hypersurface of degree $l$, $\Hilb(F)$ is naturally
isomorphic to $\P(H^0(\P^r,\O(l)))$.

\begin{lem} \label{splitting}
Let $X\subset \P^r=\P$ be a smooth surface,
$F\subset
\P$ be a  smooth hypersurface of degree $l$ transversal to $X$, and let
$C=X\cap
F$. Then there is a natural isomorphism
$N_{C/\P}|_C\iso N_{X/\P}|_C\oplus \O_C(l)$.
\end{lem}
\def\largo{\phantom{.}\\}
\begin{Pf} We have a natural diagram: $$
\begin{array}{ccccccccc}
&&&&  0  &&&&\\
&&&&  \downarrow  &&&&\\
\largo
&\phantom{\sum_I^J}&&&  N_{C/F}  &&&&\\
\largo
&\phantom{\sum_I^J}&&&  \downarrow  &\searrow&&&\\
\largo
0&\longra
&N_{C/X}&\longra&\phantom{\sum_I^J}N_{C/\P}\phantom{\sum_I^J}&\longra&N_{X/\P}|_{C}&\longra&0\\
\largo
&\phantom{\sum_I^J}&&\searrow& \downarrow  &&&&\\
\largo
&\phantom{\sum_I^J}&&&  N_{F/\P}|_{C}  &&&&\\
\largo
&&&&  \downarrow  &&&&\\
&&&&  0  &&&&
\end{array}$$
However $N_{F/\P}|_{C}=N_{C/X}=\O_C(l)$, and it is easy
to check that the morphism $N_{C/X}\to N_{F/\P}|_{C}$ in the diagram is
the identity. This implies that the natural map  $N_{C/F}\to N_{X/\P}|_{C}$ is
also an isomorphism, hence the claimed splitting.
\end{Pf}
\begin{prop}\label{smooth}
Let $X\subset \P^r=\P$ be a smooth,
regular, projectively normal surface. Let $F$ be a
smooth hypersurface of degree $l$ in $\P$ cutting
$X$ transversally along a curve $C$, and let $U\subset \Hilb(X)\times
\Hilb(F)$ be the open set of pairs $(\barx,\barf)$ such that $\barx$ and
$\barf$ are smooth and transversal and $\barx$ is projectively normal.
If $l>>0$, then for every $(\barx,\barf)\in U$ the map
$H(\barx)\times H(\barf)\to H(\barc)$ induced by intersection is smooth, where
$\barc=\barx\cap\barf$.
\end{prop}
\begin{Pf} The germ of the Hilbert scheme $H(Z)$ represents the functor of
embedded deformations of $Z$ in $\P$; when $Z$ is smooth, this functor has
tangent (resp.~obstruction) space $H^0(Z,N_{Z/\P})$ (resp.\
$H^1(Z,N_{Z/\P})$). Let $(\barx,\barf)\in U$, and $\barc=\barx\cap \barf$. The
map $H(\barx)\times H(\barf)\to H(\barc)$
induces natural maps on tangent and obstruction spaces; to prove the
required smoothness it is enough to prove that the induced maps are
surjective on tangent spaces and injective on obstruction spaces.
Note that $H^i(\barc,N_{\barc/\P})=
H^i(\barc,N_{\barx/\P}|_{\barc})\oplus H^i(\barc,N_{\barc/\P}|_{\barc})$ by
lemma \ref{splitting}. Via this isomorphism, the maps we are interested in are
induced by the long exact sequences associated to:$$
\begin{array}{ccccccccc}
0&\to &N_{\barx/\P}\otimes\I_{\barc\subset
\barx}&\longrightarrow&N_{\barx/\P}& \longrightarrow
N_{\barx/\P}|_{\barc}&\to&0\\ \phantom{1}\\
0&\to &N_{\barf/\P}\otimes\I_{\barc\subset
\barf}&\longrightarrow&N_{\barf/\P}& \longrightarrow
N_{\barf/\P}|_{\barc}&\to&0. \end{array}$$
Therefore it is enough to prove that
$$H^1(\barx,N_{\barx/\P}\otimes\I_{\barc\subset \barx})=0$$ and that
$H^0(\barf,N_{\barf/\P})\to  H^0(C,N_{\barf/\P}|_{\barc})$ is surjective
(remark that $N_{\barf/\P}=\O_\barf(l)$, hence $H^1(\barf,N_{\barf/\P})=0$ by
Kodaira vanishing).

For the claimed surjectivity, note that there is a commutative diagram $$
\begin{array}{ccc}
H^0(\P,\O_{\P}(l))&\longrightarrow &H^0(\barx,\O_{\barx}(l))\cr
\downarrow&\phantom{{2 choose 3}}&\downarrow\cr
H^0(\barf,\O_{\barf}(l))&\longrightarrow &H^0(\barc,\O_{\barc}(l))
\end{array}
$$

As $\barx$ is projectively normal, the upper horizontal arrow is onto, and as
$\barx$ is regular, the right vertical arrow is onto. Hence the lower
horizontal arrow is also onto.

To prove the vanishing, as $\I_{\barc\subset \barx}=\O_\barx(-l)$, it is enough
to prove that
$H^1(\barx,N_{\barx/\P}(-l))=0$  if $l$ is sufficiently large.  For any given
$\barx$, this follows immediately from the definition of ampleness; on the
other hand it is easy to prove (by a standard semicontinuity argument) that in
fact an $l_0$ can be found such that the claimed vanishing holds for all $l\ge
l_0$ and for all $\barx\in \Hilb(X)$.
\end{Pf}

\begin{prop}\label{inj}
Let $X\subset \P^r=\P$ be a smooth surface
and let $F$ be a smooth hypersurface of degree $l$ meeting $X$ transversally
in a smooth curve $C$. Let $U\subset \Hilb(X) \times \Hilb(F)$ be the open set
of pairs $(\barx,\barf)$ such that $\barx \cap \barf$ is a smooth curve. If
$l>>0$, then each fibre of the map $U\to \Hilb(C)$ given
by $(\barx,\barf) \mapsto \barx\cap \barf$ is contained in a fibre of the
projection $U\to\Hilb(X)$. In other words, each curve contained in the image
of $U$ in
$\Hilb(C)$ lies on exactly one surface in $\Hilb(X)$.
\end{prop}
\begin{Pf}
Let $\X\to
\Hilb(X)$ be the  universal family. Inside the product $\X \times \X$,
consider the diagonal subvariety ${\cal I}$ consisting of the pairs $(x,x)$.
Let $W$ be the locus of $\Hilb(X)\times \Hilb(X)$ over which the map ${\cal
I}\to \Hilb(X)\times \Hilb(X)$ has one-dimensional fibres. One may choose a
stratification $\{W_j\}$ of $W$ such that each of the restricted
families $\I_j\to W_j$ is flat. Thus, one has induced maps from $W_j$ to the
Hilbert scheme of one-dimensional subschemes of $\P$. Since the union of the
images of the $W_j$'s is contained in a finite number of components of the
Hibert scheme, the degree of the curves contained in the intersection of two
distinct surfaces of $\Hilb(X)$ is bounded by an integer $l_0$. Therefore it
is enough to choose $l>l_0$.
\end{Pf}

\begin{thm}\label{mainthm2}
 Let $n\ge 3$ be an integer. Then there exist
infinitely many integers $r$, and for each of them a smooth, irreducible curve
$C_r\subset \P^r$ such that $C_r$ lies exactly on $n$ components of the
Hilbert scheme of $\P^r$.
\end{thm}
\begin{Pf}
By corollary \ref{HilbX}, for infinitely many values of $r$ one can construct
a regular surface $X$ of general type, embedded in $\P^r$ by a complete
$m$-canonical system, such that $X$ lies
on exactly
$n$ components of the Hilbert scheme of $\P^r$, having pairwise different
dimensions. By \cite{And}, $X$ is projectively normal in $\P^r$ if
$m>>0$: in fact, by the theorem on page 362 together with the fact that if
$K_X$ is ample then $5K_X$ is very ample, it is enough to assume
$m\ge 11$. Choose an integer
$l$ with
$l>>0$, such that propositions
\ref{smooth} and \ref{inj} hold for $l$.

Let $F$ be a smooth hypersurface of degree $l$ meeting $X$ transversally.
Let $U\subset \Hilb(X)\times \Hilb(F)$ be the locus of pairs $(\barx,\barf)$
where both are smooth and meeting transversally, and $\barx$ is projectively
normal. $U$ is the union of
$n$ irreducible components of pairwise different dimensions, each of them being
the inverse image of an irreducible component of $\Hilb(X)$.

Let now $C=C_r$ be the intersection of $X$ and $F$, and consider the natural
map $U\to \Hilb(C)$ given by $(\barx,\barf)\mapsto \barx\cap \barf$. By
proposition
\ref{smooth} this morphism is dominant and smooth on its image $V$.  By
\ref{inj} there is an induced morphism $V\to \Hilb(X)$, which is also
dominant and smooth on its image. Therefore there is a natural
bijection between the irreducible components of $\Hilb(X)$ and those of
$\Hilb(C)$.
\end{Pf}


\bigskip
$$
\begin{array}{lll}
\hbox{\rm Dipartimento di Matematica}
&\phantom{aa}
&\hbox{\rm Dipartimento di Matematica e Informatica}
\\
\hbox{\rm Via Sommarive, 7}&&
\hbox{\rm Via Zanon, 6}
\\
\hbox{\rm I 38050 Povo - Italy}&&
\hbox{\rm I 33100 Udine - Italy}
\\
\hbox{\rm fantechi@itnvax.science.unitn.it}&&
\hbox{\rm pardini@ten.dimi.uniud.it}
\end{array}
$$

\enddocument